\documentclass{amsart}
\usepackage{amssymb}
\usepackage{amsmath}
\usepackage{amscd}
\usepackage{graphics}
\usepackage{latexsym}
\usepackage{amsrefs}
\usepackage{hyperref}

\newtheorem{theorem}{Theorem}[section]
\newtheorem{lemma}[theorem]{Lemma}

\newtheorem{remark}[theorem]{Remark}
\newtheorem{example}[theorem]{Example}

\title{Symplectic Geography Problem in Dimension Six}
\author{Ahmet Beyaz}
\address{Department of Mathematics, Middle East Technical University, Ankara 06800 Turkey}
\email{beyaz@metu.edu.tr}

\begin{document}
\begin{abstract}
In this note, the geography problem in dimension four is reviewed and then its extension to dimension six for the symplectic case is explained. Finally some examples in dimension six are provided.
\end{abstract}
\maketitle

\section{Introduction} \label{introduction}

An almost complex structure $J$ on an even dimensional smooth manifold $X^{2n}$ is a map from the tangent space $TM$ to itself such that restriction to each fiber has the property $J^2$ is negative of the identity map. If a manifold $X$ has an almost complex structure $J$ then the couple $(X,J)$ is called an almost complex manifold. We are going to drop $J$ and write $X$ for such an almost complex manifold $(X,J)$. Complex manifolds are examples of almost complex manifolds and not all almost complex manifolds are complex. An almost complex manifold $X$ has special characteristic classes, namely the Chern classes $c_i(X)$, which are elements in the even cohomology group $H^{2i}(X,\mathbb{Z})$. The zeroth Chern class of $X$ is $1\in H^0(X,\mathbb{Z})$. The total Chern class of $X$ is $c(X)=1+c_1(X)+\dots+c_n(X)$. If $X$ and $Y$ are almost complex manifolds then $c(X\times Y)$ is $c(X)c(Y)$. Here product of individual Chern classes is the cup product of cohomology classes.

If $n$ is equal to the sum of integers $r_1,\dots, r_j$, then the corresponding Chern number of $X$ is evaluation $c_{r_1}\cup\dots\cup c_{r_j}$ on the fundamental class $[X]\in H_{2n}(X,\mathbb{Z})$. The Chern numbers are integers. The top Chern number $c_n$ is equal to the Euler characteristic of the manifold. In dimension four there are only two Chern numbers which are $c_1^2$ and $c_2$. The geography problem in dimension four has two parts. First, it is the search of smooth (or symplectic) $4$-manifolds which admit a given couple $(a,b)$ of integers as the Chern numbers. The second aim is to find the number of manifolds which admit the couple $(a,b)$ as the Chern numbers. The latter is also known as the Botany problem (\cite{Fintushel2009}).

For an almost complex $4$-manifold $X$ we have the equalities $c_1^2(X)=3\sigma(X)+2e(X)$ and $c_2(X)=e(X)$ where $\sigma(X)$ is the signature of $X$ and $e(X)$ is the Euler characteristic of $X$. An almost complex $4$-manifold has another invariant which is called the holomorphic Euler characteristic $\chi_h(X)$. It turns out that $\chi_h(X)$ is  equal to $(\sigma(X)+e(X))/4$ which is completely topological. So the invariant $\chi_h(X)$ is extended to the smooth category by this equality.

Figure 1 is outlining the geography of smooth (or symplectic) $4$-manifolds which is filled by manifolds which are constructed using some state of the art techniques. More information can be found in \cite{Fintushel2009}.

The Chern numbers of an almost complex $6$-manifold $M$ are $c_3$, $c_1^3$ and $c_1c_2$. Note that $c_3(M)$ is equal to the Euler characteristic $e(M)$ of $M$. The complex geography problem in dimension six (complex dimension three) is explored by Hunt in \cite{Hunt1989}. There are some studies on almost complex geography like \cite{Okonek1995} and \cite{Lebrun1999}. We are going to focus on symplectic geography, i. e. distribution of the Chern numbers of an compatible almost complex structure on a symplectic $6$-manifold. Next theorem by Halic shows that the restrictions on the Chern numbers of symplectic $6$-manifolds are mild. 

{\theorem (Halic \cite{Halic1999}) The Chern numbers of a symplectic $6$-manifold $M$ satisfy the conditions given as:
$$c_3(M)\equiv c_1^3(M)\equiv 0\,\,(\text{mod}\,\, 2)\quad\text{and}\quad c_1c_2(M)\equiv 0\,\,(\text{mod}\,\, 24)$$
Conversely, for a triple $(2r,2s,24t)$, where $r$, $s$ and $t$ are integers, there is a simply connected symplectic $6$-manifold $M$ with property $(c_3(M),c_1^2(M),c_1c_2(M)$ is equal to $(2r,2s,24t)$.}\\

\begin{center}
\setlength{\unitlength}{1in}
\begin{picture}(4.6,4.7)
\put(.2,.5){\vector(0,1){4}}
\put(.2,.5){\vector(1,0){4.5}}
\put(.2,.5){\line(1,3){1.15}}
\put(.2,.5){\line(1,2){1.7}}
\put(.8,.5){\line(3,2){3.4}}
\put(.8,.5){\line(3,1){3.4}}
\put(0,4.3){$c_1^2$}
\put(4.6,.35){$\chi_h$}
\put(.3,.25){Elliptic Surfaces $E(n)$  ($(\chi_h,c_1^2)=(n,0)$)}
\put(3.1,.3){$ c_1^2 < 0$ (unknown)}
\put(2.0,2.7){$2\chi_h - 6 \le c_1^2 \le 9\chi_h$}
\put(1.5,2.35){surfaces of general type}
\put(1.1,4.1){$c_1^2=9\chi_h$}
\put(.3,3.1){$c_1^2>9\chi_h$}
\put(.2,2.9){(unknown)}
\put(1.8,4.1){$c_1^2=8\chi_h$}
\put(1.82,3.95){$\sigma =0$}
\put(3.9,2.8){$c_1^2=2\chi_h - 6$}
\put(3.9,1.41){$c_1^2=\chi_h - 3$}
\put(1.15,3.1){$\sigma >0$}
\put(2.25,3.1){$\sigma <0$}
\put(3.10,1.9){$\chi_h - 3 \le c_1^2 \le 2\chi_h -6$}
\put(2.65,1.65){symplectic with one basic class}
\put(3.0,1.05){$0 \le c_1^2 \le \chi_h-3$}
\put(2.0,.80){symplectic with $(\chi_h-c_1^2-2)$ basic classes}
\multiput(.170,.475)(.2,0){23}{$\bullet$}
\put (2.5, 0){Figure 1}
\end{picture}
\end{center}
\vspace{0.5cm}

Halic constructed the manifolds for the second part of the theorem by using an operation which is called generalized symplectic fiber sum. This operation was first introduced by Gompf (\cite{Gompf1995}). In six dimensional case let $M_1$ and $M_2$ be symplectic $6$-manifolds. If $X_1$ and $X_2$ are codimension two symplectic submanifolds of  $M_1$ and $M_2$ with trivial normal bundles, respectively, then the result of the fiber summing $M_1$ and $M_2$ through $X_1$ and $X_2$ is a symplectic $6$-manifold $M$. The operation is done by removing a neighborhood of $X_i$ in $M_ i$, which is symplectomorphic to $X_i\times D^2$ (\cite{Weinstein1971}). Then $M_1-nbd(X_1)$ and $M_2-nbd(X_2)$ with a symplectomorphism of $X_1\times (D^2-0)$ and $X_2\times (D^2-0)$ turning $D^2-0$ inside out. Details can be found in \cite{Gompf1995}. 

In this note we are going to apply this construction as did Halic and get some manifolds on the geography picture. Our manifolds are do not add anything new to filling the geography picture since it is already filled, but their smooth types are potentially different from those constructed by Halic. Some of these manifolds symplectic Calabi-Yau manifolds (See Remark~\ref{CY}). Recently, similar manifolds are constructed by Akhmedov in \cite{Akhmedov2011}. 

\section{Construction and Chern Numbers} \label{construction}

Let $X_1$ and $X_2$ be closed, simply connected, smooth $4$-manifolds with Lefschetz fibration structures on them. The base manifolds of the Lefschetz fibrations are $S^2$. Further information on Lefschetz fibrations may be found in \cite{Gompf1999}. Let us assume that the generic fiber of the Lefschetz fibration on the respective $X_i$ is an embedded surface $\Sigma_i$ which has genus $g_i$. Assume that each singular fiber has a unique double point (fishtail) singularity. Let $n_i$ be the number of singular fibers on each Lefschetz fibration. The Euler characteristic $e(X_i)$ of such a $4$-manifold is $2(2-2g_i)+n_i$. In particular $e(X_i)$ is equal to $c_2(X_i)$. $X_i$ is simply connected means that $n_i$ is greater than $2g_i$.

Consider the product $6$-manifolds $X_1\times\Sigma_2$ and $X_2\times\Sigma_1$ with the product symplectic structures on them. Let us apply the symplectic fiber sum operation (\cite{Gompf1995}, \cite{Halic1999}) to these $6$-manifolds over $\Sigma_1\times\Sigma_2$. As a result we get a closed, symplectic $6$-manifold which we are going to denote by $M$. 

As a symplectic manifold, $M$ has a contractible family of almost complex structures which are compatible with the symplectic structure on it. With the formula of the Chern numbers for symplectic fiber sum (\cite{Halic1999}), calculations for the Chern numbers $c_1^3$, $c_1c_2$ and $c_3$ of a generic almost complex structure on the product manifold $M$ are given as follows:

\begin{displaymath}
\begin{array}{lll} c_3(M)&=&c_3(X_1\times\Sigma_2)+c_3(X_2\times\Sigma_1)-2c_2(\Sigma_1\times\Sigma_2)\\
&=&c_2(X_1)c_1(\Sigma_2)+c_2(X_2)c_1(\Sigma_1)-2c_1(\Sigma_1)c_1(\Sigma_2)\\
&=&c_2(X_1)(2-2g_2)+c_2(X_2)(2-2g_1)-2(2-2g_1)(2-2g_2)\\
&=&(2(2-2g_1)+n_1)(2-2g_2)+(2(2-2g_1)+n_2)(2-2g_2)\\
&&-2(2-2g_1)(2-2g_2)\\
&=&8(1-g_1)(1-g_2)+2n_1(1-g_2)+2n_2(1-g_1)
\end{array}
\end{displaymath}

Considering $c_2(X_i)=e(X_i)=12\chi_h(X_i)-c_1^2(X_i)$, we can write $c_3$ in terms of $\chi_h$ and $c_1^2$. 

\begin{displaymath}
\begin{array}{lll} c_3(M)&=&2(12\chi_h(X_1)-c_1^2(X_1))(1-g_2)\\
&&+2(12\chi_h(X_2)-c_1^2(X_2))(1-g_1)-8(1-g_1)(1-g_2)
\end{array}
\end{displaymath}

Calculation of $c_1^3$:

\begin{displaymath}
\begin{array}{lll} c_1^3(M)&=&c_1^3(X_1\times\Sigma_2)+c_1^3(X_2\times\Sigma_1)-6c_1^2(\Sigma_1\times\Sigma_2)\\
&=&(c_1(X_1)+c_1(\Sigma_2))^3+(c_1(X_2)+c_1(\Sigma_1))^3-12c_1(\Sigma_1)c_1(\Sigma_2)\\
&=&c_1^3(X_1)+3c_1^2(X_1)c_1(\Sigma_2)+3c_1(X_1)c_1^2(\Sigma_2)+c_1^3(\Sigma_2)\\
&&+c_1^3(X_2)+3c_1^2(X_2)c_1(\Sigma_1)+3c_1(X_2)c_1^2(\Sigma_1)+c_1^3(\Sigma_1)\\
&&-12c_1(\Sigma_1)c_1(\Sigma_2)\\
&=&6(1-g_2)c_1^2(X_1)+6(1-g_1)c_1^2(X_2)-48(1-g_1)(1-g_2)
\end{array}
\end{displaymath}

Calculation of $c_1c_2$:

\begin{displaymath}
\begin{array}{lll} c_1c_2(M)&=&c_1c_2(X_1\times\Sigma_2)+c_1c_2(X_2\times\Sigma_1)-2c_1^2(\Sigma_1\times\Sigma_2)-2c_2(\Sigma_1\times\Sigma_2)\\
&=&(c_1(X_1)+c_1(\Sigma_2))(c_2(X_1)+c_1(X_1)c_1(\Sigma_2))\\
&&+(c_1(X_2)+c_1(\Sigma_1))(c_2(X_2)+c_1(X_2)c_1(\Sigma_1))\\
&&-2(c_1(\Sigma_1)+c_1(\Sigma_2))^2-2c_1(\Sigma_1)c_1(\Sigma_2)\\
&=&c_1^2(X_1)c_1(\Sigma_2)+c_2(X_1)c_1(\Sigma_2)+c_1^2(X_2)c_1(\Sigma_1)+c_2(X_2)c_1(\Sigma_1)\\
&&-4c_1(\Sigma_1)c_1(\Sigma_2)-2c_1(\Sigma_1)c_1(\Sigma_2)\\
&=&(2-2g_2)(c_1^2(X_1)+c_2(X_1))+(2-2g_1)(c_1^2(X_2)+c_2(X_2))\\
&&-6(2-2g_1)(2-2g_2)\\
&=&(2-2g_2)(3\sigma(X_1)+3e(X_1))+(2-2g_1)(3\sigma(X_2+3e(X_2))\\
&&-6(2-2g_1)(2-2g_2)
\end{array}
\end{displaymath}

If we summarize, the Chern numbers of $M$ are as follows.

{\lemma \label{Chern numbers} The Chern numbers of the symplectic $6$-manifold $M$ constructed by fiber summing $X_1\times\Sigma_2$ and $X_2\times\Sigma_1$ along $\Sigma_1\times\Sigma_2$ are:

\begin{displaymath}
\begin{array}{lll} c_3(M)&=&2(12\chi_h(X_1)-c_1^2(X_1))(1-g_2)\\
&&+2(12\chi_h(X_2)-c_1^2(X_2))(1-g_1)-8(1-g_1)(1-g_2)\\
c_1^3(M)&=&6(1-g_2)c_1^2(X_1)+6(1-g_1)c_1^2(X_2)-48(1-g_1)(1-g_2)\\
c_1c_2(M)&=&24(1-g_2)\chi_h(X_1)+24(1-g_1)\chi_h(X_2)-24(1-g_1)(1-g_2)
\end{array}
\end{displaymath}

where $\chi_h(X_i)$ is the holomorphic Euler characteristic of $X_i$.}\\

\section{Examples} \label{examples}

In this section we give two families of symplectic $6$-manifolds which are constructed using the fiber summing as above. 

{\example \label{example1} \textup{Let $X_1$ be the elliptic surface $E(m)$ with the Lefschetz fibration structure whose generic fibers are tori. $\chi_h(X_1)$ is $m$ and $c_1^2(X_1)$ is $0$. Let $X_2$ be the ruled surface $S^2\times S^2$ and consider the sphere fibration structure on it with sphere sections of self-intersection zero. We may think of it as a Lefschetz fibration with no singular fibers. Genus $g_2$ is zero. $\chi_h(X_2)$ is $1$ and $c_1^2(X_2)$ is $8$. Apply the construction explained in Section \ref{construction} to manifolds $X_1$, $X_2$. Then the Chern numbers of the constructed symplectic manifold $M$ are $c_3(M)=c_1c_2(M)=24m$ and $c_1^3(M)=0$.}}

{\example \label{example2} \textup{Let $Y_1$ be the elliptic surface $E(m)$ ($m\geq 1$). Let $Y_2$ be the manifold $E(k)_K$ ($k\geq 1$) which is obtained by applying the knot surgery of Fintushel-Stern (\cite{Fintushel1998a}) to the elliptic surface $E(k)$ over the fibered knot $K$. For a fixed $k$, $E(k)$ and $E(k)_K$ are homeomorphic to each other which means that $\chi_h(Y_2)$ is $k$ and $c_1^2(Y_2)$ is zero. $Y_2$ has an Lefschetz fibration with generic fiber genus $2g+k-1$, where $g\geq 0$ is the genus of the knot (\cite{Fintushel2004c}). If we apply the construction to manifolds $Y_1$, $Y_2$, then the Chern numbers of the constructed symplectic manifold $M$ are $c_3(M)=c_1c_2(M)=24m(2-2g-k)$ and $c_1^3(M)=0$.}}

{\remark \label{CY} All examples given above are simply connected. Throughout these families, although $c_1^3$ is zero, the only manifold with $c_1$ equal to zero is when we choose $m=k=2$ and $g=0$ in Example~\ref{example2}. See \cite{Akhmedov2011} for further discussion on the symplectic Calabi-Yau manifolds.}

{\remark Many similar examples may be produced by applying the procedure above. Still there will be triples $(2r,2s,24t)$ which we can not realize as Chern numbers due to the fact that $c_1^3$ is divisible by six for the product manifold of this operation. This is behind the result of Halic in \cite{Halic1999}. However, potentially all of manifolds constructed in Example~\ref{example2} are diffeomorphic for a fixed $k$. It may be interesting to see if they are diffeomorphic.}

\section*{Acknowledgements}

The author would like to thank Ron Stern for suggesting this construction.


\end{document}